\documentclass[a4paper,12pt,reqno]{amsart}
\usepackage{amsfonts}
\usepackage{amsmath}
\usepackage{amssymb}
\usepackage[a4paper]{geometry}
\usepackage{mathrsfs}
\usepackage{csquotes}
\usepackage[colorlinks]{hyperref}
\renewcommand\eqref[1]{(\ref{#1})} 
\numberwithin{equation}{section}
\theoremstyle{plain}
\newtheorem{thm}{Theorem}[section]

\theoremstyle{definition}

\newtheorem{ex}[thm]{Example}

\begin{document}

   \title[Isoperimetric inequalities for some integral operators]
   {Isoperimetric inequalities for some integral operators arising in potential theory}
 
\author[M. Ruzhansky]{Michael Ruzhansky}
\address{
	Michael Ruzhansky:
	\endgraf
	Department of Mathematics
	\endgraf
	Imperial College London
	\endgraf
	180 Queen's Gate, London SW7 2AZ
	\endgraf
	United Kingdom
	\endgraf
	{\it E-mail address} {\rm m.ruzhansky@imperial.ac.uk}
}

\author[D. Suragan]{Durvudkhan Suragan}
\address{
  Durvudkhan Suragan:
  \endgraf
  Institute of Mathematics and Mathematical Modelling
  \endgraf
  125 Pushkin str.
  \endgraf
  050010 Almaty
  \endgraf
  Kazakhstan
  \endgraf
  {\it E-mail address} {\rm suragan@math.kz}
  }

\thanks{The authors were supported in parts by the EPSRC
 grant EP/K039407/1 and by the Leverhulme Grant RPG-2014-02. No new data was collected or
generated during the course of research.}


     \begin{abstract}
    In this paper we review our previous isoperimetric 
    results for the logarithmic potential and Newton potential operators. The main reason why the results are useful, beyond the intrinsic interest of geometric extremum problems, is that they produce \emph{a priori} bounds for spectral invariants of operators on arbitrary domains. We demonstrate these in explicit examples. 
     \end{abstract}
     \maketitle

\section{Introduction}

In a bounded domain of the Euclidean space $\Omega\subset \mathbb{R}^{d},\,\, d\geq 2,$ it is well known that the
solution to the Laplacian equation
\begin{equation}
-\Delta u(x)=f(x), \,\,\,\ x\in\Omega,\label{t15}
\end{equation}
is given by the Newton potential formula (or the logarithmic potential formula when $d=2$)
\begin{equation}
u(x)=\int_{\Omega}\varepsilon_{d}(|x-y|)f(y)dy,\,\, x\in\Omega,
\label{t14}
\end{equation}
for suitable functions $f$ with $supp f\subset\Omega$.
Here
\begin{equation}
\varepsilon_{d}(|x-y|)=\left\{
\begin{array}{ll}
\frac{1}{2\pi}\ln\frac{1}{|x-y|}, \,\,d=2, \\
\frac{1}{(d-2)s_{d}}\frac{1}{|x-y|^{d-2}},\,\,d\geq 3,\\
\end{array}
\right.
\label{tEQ:fs}
\end{equation}
is the fundamental solution to $-\Delta$ and $s_{d}=\frac{2\pi^{\frac{d}{2}}}{\Gamma(\frac{d}{2})}$ is the surface area of the unit sphere in
$\mathbb{R} ^{d}$.

An interesting  question having several important applications is what boundary conditions can be put on $u$ on
the (Lipschitz) boundary $\partial\Omega$ so that equation \eqref{t15} complemented by this
boundary condition would have the solution in $\Omega$ still given by the same formula \eqref{t14},
with the same kernel $\varepsilon_d$ given by \eqref{tEQ:fs}.
It turns out that the answer to this question is the integral boundary condition
\begin{equation}
-\frac{1}{2}u(x)+\int_{\partial\Omega}\frac{\partial\varepsilon_{d}(|x-y|)}{\partial n_{y}}u(y)d S_{y}-
\int_{\partial\Omega}\varepsilon_{d}(|x-y|)\frac{\partial u(y)}{\partial n_{y}}d S_{y}=0,\,\,
x\in\partial\Omega,
\label{t16}
\end{equation}
where $\frac{\partial}{\partial n_{y}}$ denotes the outer normal
derivative at a point $y$ on $\partial\Omega$.
A converse question to the one above would be to determine the trace of the Newton potential
\eqref{t14} on the boundary surface $\partial\Omega$, and one can use the potential theory to show that it
has to be given by \eqref{t16}.

The boundary condition \eqref{t16} appeared in M. Kac's work \cite{Kac1} where
he called it
``the principle of not feeling the boundary''. This was further expanded in Kac's book \cite{Kac2}
with several further applications to the spectral theory and the asymptotics of the Weyl's eigenvalue
counting function.
Independently in \cite{KS1} T.Sh. Kal'menov and the second author proved the existence of the boundary condition \eqref{t16} and as byproduct the eigenvalues and eigenfunctions of the Newton potential
\eqref{t14} were calculated in the 2-disk and in the 3-ball.
In general, the boundary value problem \eqref{t15}-\eqref{t16} has various interesting properties and applications
(see, for example, \cite{Kac1,Kac2}, \cite{BS}, \cite{Sa},\cite{KS1}, \cite{Ruzhansky-Suragan:PAMS} and \cite{Ruzhansky-Suragan:Layers}).
The boundary value problem \eqref{t15}-\eqref{t16} can also be generalised for higher degrees of the
Laplacian, see \cite{KS2, KS3}.
In the present paper we consider spectral problems of inverse operators to the nonlocal Laplacian \eqref{t15}-\eqref{t16}, namely the logarithmic potential operator on $L^{2}(\Omega)$ defined by
\begin{equation}
\mathcal{L}_{\Omega}f(x):= \int_{\Omega}\frac{1}{2\pi}\ln\frac{1}{|x-y|}f(y)dy,\quad f\in L^{2}(\Omega),\quad \Omega\subset \mathbb R^{2}, \label{3}
\end{equation}
and the Newton potential operator on $L^{2}(\Omega)$ defined by
\begin{equation}
\mathcal{N}_{\Omega}f(x):= \int_{\Omega} \frac{1}{(d-2)s_{d}}\frac{1}{|x-y|^{d-2}}f(y)dy,\quad f\in L^{2}(\Omega),\quad \Omega\subset \mathbb R^{d},\,d\geq 3.\label{3}
\end{equation}
Spectral properties of the logarithmic and the Newton potential operator have been considered in many papers (see, e.g. \cite{AKL}, \cite{AK}, \cite{BS}, \cite{Kac3}, \cite{KS3}, \cite{Tr67}, \cite{Tr69} and \cite{Zo}). In this paper we are interested in isoperimetric inequalities of these operators, that is also, in isoperimetric inequalities of the nonlocal Laplacian \eqref{t15}-\eqref{t16}.   A recent general review of isoperimetric inequalities for the Dirichlet, Neumann and other Laplacians was made by Benguria, Linde and Loewe in \cite{Ben}. In addition to \cite{Ben}, we refer G.~P{\'o}lya and G.~Szeg{\"o} \cite{Po1}, Bandle \cite{Ban80} and  Henrot \cite{He} for historic remarks on isoperimetric inequalities, namely the Rayleigh-Faber-Krahn inequality and the Luttinger inequality.

We review an analogue of the Luttinger inequality for the Newton potential operator $\mathcal{N}_{\Omega}$ and provide related explicit examples. It is a particular case of our previous result with G. Rozenblum in \cite{Ruzhansky-Suragan:Newton} for the Newton potential (see also \cite{Ruzhansky-Suragan:UMN},   \cite{Ruzhansky-Suragan:BMS} and \cite{KS16} for a non-self adjoint operators).
In Section \ref{sec:3} we present:
\begin{itemize}
	\item Luttinger type inequality for $\mathcal{N}_{\Omega}$: The $d$-ball is a maximizer of the Schatten $p$-norm of the Newton potential operator among all domains of a given measure in $\mathbb R^{d},\,d\geq 3$, for all integer $\frac{d}{2}<p<\infty$.
\end{itemize}
In Section \ref{sec:2}, we review the following facts for the logarithmic potential from \cite{Ruzhansky-Suragan:arxiv}:
\begin{itemize}
	\item Rayleigh-Faber-Krahn inequality: The disc is a minimizer of the characteristic number of the logarithmic potential $\mathcal{L}_\Omega$ with the smallest modulus among all domains of a given measure.
	\item P{\'o}lya inequality: The equilateral triangle is a minimizer of the first characteristic number of the logarithmic potential $\mathcal{L}_\Omega$ with the smallest modulus among all triangles of a given area.
	\item Luttinger type inequality for $\mathcal{L}_{\Omega}$: The disc is a maximizer of the Schatten $p$-norm of the logarithmic potential operator among all domains of a given measure in $\mathbb R^{2}$, for all integer $2\leq p<\infty$.
	\item Luttinger type inequality for $\mathcal{L}_{\Omega}$ in triangles: The equilateral triangle is a maximizer of the Schatten $p$-norm of the logarithmic potential operator among all triangles of a given area in $\mathbb R^{2}$, for all integer $2\leq p<\infty$.
\end{itemize}

\section{Isoperimetric inequalities for $\mathcal{L}_\Omega$ and examples}
\label{sec:2}
In this section we review our results for the logarithmic potential from \cite{Ruzhansky-Suragan:arxiv}
and for the Newton potential \cite{Ruzhansky-Suragan:Newton}.
Let $\Omega\subset \mathbb R^{2}$ be an open bounded set. We consider the logarithmic potential operator on $L^{2}(\Omega)$ defined by
\begin{equation}
\mathcal{L}_{\Omega}f(x):= \int_{\Omega}\frac{1}{2\pi}\ln\frac{1}{|x-y|}f(y)dy,\quad f\in L^{2}(\Omega),\label{3}
\end{equation}
where $\ln$ is the natural logarithm and $|x-y|$ is the standard Euclidean
distance between $x$ and $y$.
Clearly, $\mathcal{L}_{\Omega}$ is compact and self-adjoint.
Therefore, all of its eigenvalues and characteristic numbers are discrete and real.
We recall that the characteristic numbers are inverses of the eigenvalues. The characteristic numbers of $\mathcal{L}_{\Omega}$ may be enumerated in ascending order of their modulus,
$$|\mu_{1}(\Omega)|\leq|\mu_{2}(\Omega)|\leq...$$
where $\mu_{i}(\Omega)$ is repeated in this series according to its multiplicity. We denote the
corresponding eigenfunctions by $u_{1}, u_{2},...,$ so that for each characteristic number $\mu_{i}$ there is a unique corresponding (normalized) eigenfunction $u_{i}$,
$$u_{i}=\mu_{i}(\Omega)\mathcal{L}_{\Omega}u_{i},\,\,\,\, i=1,2,....$$
\medskip
It is known, see for example \cite{KS1},
that the equation
$$
u(x)=\mathcal{L}_{\Omega} f(x)=\int_{\Omega}\frac{1}{2\pi}\ln\frac{1}{|x-y|}f(y)dy
$$
is equivalent to the equation
\begin{equation}\label{15}
-\Delta u(x)=f(x), \,\,\,\ x\in\Omega,
\end{equation}
with the nonlocal integral boundary condition
\begin{equation}
-\frac{1}{2}u(x)+\int_{\partial\Omega}\frac{\partial}{\partial n_{y}}\frac{1}{2\pi}\ln\frac{1}{|x-y|}u(y)d S_{y}-
\int_{\partial\Omega}\frac{1}{2\pi}\ln\frac{1}{|x-y|}\frac{\partial u(y)}{\partial n_{y}}d S_{y}=0,\,\,
x\in\partial\Omega,
\label{16}
\end{equation}
where $\frac{\partial}{\partial n_{y}}$ denotes the outer normal
derivative at a point $y$ on the boundary $\partial\Omega$, which is assumed piecewise $C^1$ here.

Let $H$ be a separable Hilbert space. By $\mathcal{S}^{\infty}(H)$
we denote the space of compact operators $P:H\rightarrow H$.
Recall that the singular numbers $\{s_{n}\}$ of $P\in \mathcal{S}^{\infty}(H)$
are the eigenvalues of the positive operator $(P^{*}P)^{1/2}$ (see \cite{GK}). The Schatten $p$-classes are defined as
$$
\mathcal{S}^{p}(H):=\{P \in \mathcal{S}^{\infty}(H): \{s_{n}\}\in \ell^{p}\},\quad 1\leq p<\infty.
$$
In $\mathcal{S}^{p}(H)$ the Schatten $p$-norm of the operator $P$ is defined by
\begin{equation}
\| P\|_{p}:=\left(\sum_{n=1}^{\infty}s_{n}^{p}\right)^{\frac{1}{p}}, \quad 1\leq p <\infty.\label{1}
\end{equation}
For $p=\infty$, we can set
$$
\| P\|_{\infty}:=\| P\|
$$
to be the operator norm of $P$ on $H$. As outlined in the introduction, we assume that $\Omega\subset \mathbb R^{2}$
is an open bounded set and we consider the logarithmic potential operator
on $L^{2}(\Omega)$ of the form
\begin{equation}
\mathcal{L}_{\Omega}f(x)=\int_{\Omega}\frac{1}{2\pi}\ln\frac{1}{|x-y|}f(y)dy,\quad f\in L^{2}(\Omega).
\label{n16}
\end{equation}
It is known that $\mathcal{L}_{\Omega}$ is a Hilbert-Schmidt operator.
By $|\Omega|$ we will denote the Lebesque measure of $\Omega$.

\begin{thm} \label{THM:main}
	Let $D$ be a disc centred at the origin.
	Then
	\begin{equation}
	\|\mathcal{L}_{\Omega}\|_{p}\leq  \|\mathcal{L}_{D}\|_{p}
	\label{3}
	\end{equation}
	for any integer $2\leq p\leq \infty$ and any bounded open domain $\Omega$ with $|\Omega|=|D|.$
\end{thm}

Let us give several examples calculating explicitly values of the right hand side of \eqref{3} for different values of $p$.
\begin{ex}
	Let $D\equiv U$ be the unit disc.
	Then by Theorem \ref{THM:main} we have
	\begin{equation}
	\|\mathcal{L}_{\Omega}\|_{p}\leq
	\|\mathcal{L}_{U}\|_{p}=
	\left(\sum_{m=1}^{\infty}\frac{3}{j_{0,m}^{2p}}+\sum_{l=1}^{\infty}\sum_{m=1}^{\infty}\frac{2}{j_{l,m}^{2p}}\right)^{\frac{1}{p}},
	\label{ex2}
	\end{equation}
	for any integer $2\leq p < \infty$ and any bounded open domain $\Omega$ with $|\Omega|=|U|.$
	Here $j_{km}$ denotes the $m^{th}$ positive zero of the Bessel function  $J_{k}$ of the first kind of order $k$.
	
	The right hand sight of the formula \eqref{ex2} can be confirmed by a direct calculation of the logarithmic potential eigenvalues in the unit disc, see Theorem 3.1 in \cite{AKL}.
\end{ex}
\begin{ex}
	Let $D\equiv U$ be the unit disc.
	Then by Theorem \ref{THM:main} we have
	\begin{equation}
	\| \mathcal{L}_\Omega\|
	\leq \| \mathcal{L}_U\| =\frac{1}{j_{01}^{2}}
	\label{ex3}
	\end{equation}
	for any bounded open domain $\Omega$ with $|\Omega|=|U|.$ Here $\|\cdot\|$ is the operator norm on the space $L^2$.
\end{ex}
From Corollary 3.2 in \cite{AKL} we calculate explicitly the operator norm in the right hand sight of
\eqref{ex3}.

In Theorem \ref{THM:main} when $p=\infty$, the following analogue of the Rayleigh-Faber-Krahn theorem for the integral operator $\mathcal{L}_\Omega$ is used.

\begin{thm} \label{THM:1} The disc $D$ is a minimizer of the characteristic number of the logarithmic potential $\mathcal{L}_\Omega$ with the smallest modulus among all domains of a given measure, i.e.
	$$0<|\mu_{1}(D)|\leq |\mu_{1}(\Omega)|$$
	for an arbitrary bounded open domain $\Omega\subset \mathbb R^{2}$ with $|\Omega|=|D|.$
\end{thm}

In Landkof \cite{Lan66} the positivity of the operator $\mathcal{L}_\Omega$ is proved in domains $\overline{\Omega}\subset U,$ where $U$ is the unit disc. In general, $\mathcal{L}_\Omega$ is not a positive operator. For any bounded open domain $\Omega$ the logarithmic potential operator $\mathcal{L}_\Omega$ can have at most one negative eigenvalue, see
Troutman \cite{Tr67} (see also Kac \cite{Kac3}).

In other words Theorem \ref{THM:1} says that the operator norm of
$\mathcal{L}_\Omega$ is maximized in a disc among all Euclidean bounded open domains of a given area.

It follows from the properties of the kernel that the Schatten $p$-norm of the operator $\mathcal{L}_\Omega$ is finite when $p>1$ see e.g. the criteria for Schatten classes in terms of the regularity of the kernel in \cite{DR}.
Our techniques do not allow us to prove Theorem \ref{THM:main} for $1< p<2$. In view of the Dirichlet Laplacian case, it seems reasonable to conjecture that the Schatten $p$-norm is still maximized on the disc also for $1< p<2$. However, In Section 3 by using different method we prove such conjecture for the Newton potential operator, see also \cite{Ruzhansky-Suragan:Newton}.

We can ask the same question of maximizing the Schatten $p$-norms in the class of polygons with a
given number $n$ of sides. We denote by $\mathcal{P}_{n}$ the class of plane polygons with $n$ edges.
We would like to identify the maximizer for Schatten $p$-norms of the logarithmic potential $\mathcal{L}_\Omega$ in $\mathcal{P}_{n}$.
According to the Dirichlet Laplacian case, it is natural to conjecture that it is the $n$-regular polygon. Currently, we have proved this only for $n=3$:

\begin{thm}\label{THM:tri}
	The equilateral triangle centred at the origin has the largest Schatten $p$-norm of the operator
	$\mathcal{L}_\Omega$ for any integer $2\leq p\leq \infty$ among all triangles of a given area.
	More precisely, if $\triangle$ is the equilateral triangle centred at the origin, we have
	\begin{equation}
	\|\mathcal{L}_{\Omega}\|_{p}\leq  \|\mathcal{L}_{\triangle}\|_{p}
	\label{3-1}
	\end{equation}
	for any integer $2\leq p\leq \infty$ and any bounded open triangle $\Omega$ with $|\Omega|=|\triangle|.$
\end{thm}

When $p=\infty$, Theorem \ref{THM:tri} implies the following analogue of the P{\'o}lya theorem \cite{Po} for the operator $\mathcal{L}_\Omega$.

\begin{thm} \label{THM:2}
	The equilateral triangle  $\triangle$ centred at the origin is a minimizer of the first characteristic number of the logarithmic potential $\mathcal{L}_\Omega$ among all triangles of a given area, i.e.
	$$0<|\mu_{1}(\triangle)|\leq |\mu_{1}(\Omega)|$$
	for any triangle $\Omega\subset \mathbb R^{2}$ with $|\Omega|=|\triangle|.$
\end{thm}

In other words Theorem \ref{THM:2} says that the operator norm of
$\mathcal{L}_\Omega$ is maximized in an equilateral triangle among all triangles of a given area.

\section{The Newton potential}
\label{sec:3}
Let $\Omega\subset  \mathbb R^{d},\,\,d\geq 3,$ be an open bounded set.
We consider the Newton potential operator $\mathcal{N}_{\Omega}:L^{2}(\Omega)\rightarrow L^{2}(\Omega)$ defined by
\begin{equation}
\mathcal{N}_{\Omega}f(x):= \int_{\Omega}\varepsilon_{d}(|x-y|)f(y)dy,\quad f\in L^{2}(\Omega),
\label{6.3}
\end{equation}
where $\varepsilon_{d}(|x-y|)=\frac{1}{(d-2)s_{d}}\frac{1}{|x-y|^{d-2}},\,\,d\geq 3.$

Since $\varepsilon_{d}$ is positive, real and symmetric function, $\mathcal{N}_{\Omega}$ is a positive self-adjoint operator.
Therefore, all of its eigenvalues and characteristic numbers are positive real numbers.
We recall that the characteristic numbers are inverses of the eigenvalues.
The characteristic numbers of $\mathcal{N}_{\Omega}$ may be enumerated in ascending order

$$0<\mu_{1}(\Omega)\leq\mu_{2}(\Omega)\leq\ldots,$$
where $\mu_{i}(\Omega)$ is repeated in this series according to its multiplicity. We denote the
corresponding eigenfunctions by $u_{1}, u_{2},...,$ so that for each characteristic number $\mu_{i}$ there is a unique corresponding (normalized) eigenfunction $u_{i}$,
$$u_{i}=\mu_{i}(\Omega)\mathcal{N}_{\Omega}u_{i},\,\,\,\, i=1,2,\ldots.$$

This spectral problem has various interesting properties and applications (see \cite{Kac1} and \cite{Sa}, for example). In particular, one can prove that in the unit ball its spectrum contains the spectrum of the corresponding Dirichlet Laplacian by using an explicit calculation (cf. \cite{KS3}).

Kac \cite{Kac1} proved that
\begin{equation}
1=\lim_{\delta\rightarrow 0} \sum_{j=1}^{\infty}\frac{1}{1+\mu_{j}\delta}u_{j}(y)\int_{\Omega}u_{j}(x)dx, \quad
y\in \Omega,\label{eq4}
\end{equation}
where $\mu_{j}, j=1,2,...,$ and $u_{j}, j=1,2,...,$ are the characteristic numbers and the corresponding normalized eigenfunctions of the Newton potential operator \eqref{6.3}, respectivelty.
The purely analytic fact \eqref{eq4} expresses that the expansion of $1$ in a series of orthonormal functions $u_{j}$ is summable to $1$ for every $y\in \Omega$. In \cite{Kac4} Kac gave asymptotic formulae for the characteristic numbers in $\mathbb R^{d}, d\geq3$. 
In this section we discuss some other pure analytic facts for the Newton potential. It should be noted that similar results are already known for the Dirichlet Laplacian. 

By using the Feynman-Kac formula and spherical rearrangement, Luttinger proved that the ball $\Omega^{*}$ is the maximizer
of the partition function of the Dirichlet Laplacian among all domains of the same volume as $\Omega^{*}$ for all positive values of time \cite{Lu}, i.e.

\begin{equation*}
Z_{\Omega}^{\mathcal{D}}(t):=
\sum_{i=1}^{\infty}\exp(-t\lambda_{i}^{\mathcal{D}}(\Omega))\leq Z_{\Omega^{*}}^{\mathcal{D}}(t):=\sum_{i=1}^{\infty}\exp(-t\lambda_{i}^{\mathcal{D}}(\Omega^{*})), \quad |\Omega|=|\Omega^{*}|,\forall t>0,\end{equation*}
where $\lambda_{i}^{\mathcal{D}}(\Omega), i=1,2,...,$ are the eigenvalues of the Dirichlet Laplacian $\Delta^{\mathcal{D}}_{\Omega}$ in $\Omega$.

The partition function and the Schatten norms are related:
\begin{equation*}\|\Delta^{\mathcal{D}}_{\Omega}\|_{p}^{p}
=\frac{1}{\Gamma(p)}\int_{0}^{\infty}t^{p-1}Z_{\Omega}^{\mathcal{D}}(t)dt,\end{equation*}
where $\Gamma$ is the gamma function.
Hence it easily follows that
\begin{equation}\|\Delta^{\mathcal{D}}_{\Omega}\|_{p}\leq
\|\Delta^{\mathcal{D}}_{\Omega^{*}}\|_{p},\,\,|\Omega|=|\Omega^{*}|, \label{eq5}\end{equation}
when $p>d/2$, $\Omega\subset \mathbb R^{d}$.
Here the Schatten $p$-norm of the Dirichlet Laplacian is defined by
$$\|\Delta^{\mathcal{D}}_{\Omega}\|_{p}:=\left(\sum_{i=1}^{\infty}
\frac{1}{[\lambda_{i}^{\mathcal{D}}]^{p}}\right)^{p},\quad d/2<p<\infty.$$
The right hand side of the inequality \eqref{eq5} gives the exact upper bound of the Schatten $p$-norm and it can be calculated explicitly.

\begin{ex} Let $U$ be the unit disk then
	\begin{equation*}\|\Delta^{\mathcal{D}}_{U}\|_{2}^{2}=0.0493....\end{equation*}
	Therefore, from \eqref{eq5} we have
	\begin{equation*}\|\Delta^{\mathcal{D}}_{\Omega}\|_{2}^{2}\leq 0.0493..., \,\,|\Omega|=|U|.\end{equation*}
	This inequality is better than the inequality conjectured in \cite{HH}
	\begin{equation}\|\Delta^{\mathcal{D}}_{\Omega}\|_{p}^{p}\leq \frac{\Gamma(p-\frac{d}{2})}{\Gamma(p)}\frac{Vol|\Omega|
		^{\frac{2p}{d}}}{(4\pi)^{\frac{d}{2}}},\,\, p>\frac{d}{2},\label{eq6}\end{equation}
	which implies that
	\begin{equation*}\|\Delta^{\mathcal{D}}_{\Omega}\|_{2}^{2}\leq 0.7853...,\end{equation*}
	when $|\Omega|=|U|$.
\end{ex}
However, it is important to note that in \eqref{eq6} $p$ is an arbitrary real number greater than $\frac{d}{2}$.

The condition $p>d/2$ in \eqref{eq5} is necessary to absolute convergence of series,
but in case $p\leq d/2$ one may use regularization process to get an absolute convergent series.

\begin{ex}[\cite{Dos11}] In $\Omega\subset \mathbb R^{2}$ the sum
	\begin{equation*}\|\Delta^{\mathcal{D}}_{\Omega}\|_{1}
	=\sum_{k=1}^{\infty}\frac{1}{\lambda_{k}^{\mathcal{D}}(\Omega)}=\infty, \,\,\Omega\subset \mathbb R^{2}.\end{equation*}
	However, using the following regularisation we find that
	if $U\equiv\Omega\subset R^{2}$ is the unit disk then
	\begin{equation*}
	\sum_{k=1}^{\infty}
	\left(\frac{1}{\lambda_{k}^{\mathcal{D}}(U)}-\frac{1}{4k}\right)=
	-0.3557....\end{equation*}
\end{ex}

As usual by $|\Omega|$ we will denote the Lebesque measure of $\Omega$.

\begin{thm} \label{THM:N}
	Let $B$ be a ball centred at the origin, $d\geq 3$.
	Then
	\begin{equation}
	\|\mathcal{N}_{\Omega}\|_{p}\leq  \|\mathcal{N}_{B}\|_{p}
	\label{3-1}
	\end{equation}
	for any integer $\frac{d}{2}< p\leq \infty$ and an arbitrary bounded open domain $\Omega$ with $|\Omega|=|B|.$
\end{thm}

Let us give some examples:
\begin{ex}
	Let $B\equiv U$ be the unit 3-ball.
	Then by Theorem \ref{THM:N} we have
	\begin{equation}
	\|\mathcal{N}_{\Omega}\|_{p}\leq
	\|\mathcal{N}_{U}\|_{p}=
	\left(\sum_{l=0}^{\infty}\sum_{m=1}^{\infty}\frac{2l+1}{j_{l-\frac{1}{2},m}^{2p}}\right)^{\frac{1}{p}},
	\label{ex2}
	\end{equation}
	for any real $2\leq p < \infty$ and any bounded open domain $\Omega$ with $|\Omega|=|U|.$
	Here $j_{km}$ denotes the $m^{th}$ positive zero of the Bessel function  $J_{k}$
	of the first kind of order $k$.
	The right hand side of the formula \eqref{ex2} can be confirmed by a direct calculation of the characteristic numbers of the Newton potential in the unit 3-ball,
	see Theorem 4.1 in \cite{AKL} (cf. \cite{KS1}).
\end{ex}

\begin{ex}
	For the Hilbert-Schmidt norm we have
	\begin{equation}
	\|\mathcal{N}_{\Omega}\|_{2}\leq
	\|\mathcal{N}_{U}\|_{2}=
	\sqrt{\frac{7}{48}},\label{ex}
	\end{equation}
	for any bounded open domain $\Omega$ with $|\Omega|=|U|,$ where $B\equiv U$ is the unit 3-ball.
	Here, when $p=2$, we have calculated the value on the right hand side
	of the inequality \eqref{ex} by using
	the polar representation. We omit the routine technical calculation.
\end{ex}
\begin{ex}
	When $p=\infty$ by Theorem \ref{THM:N} we have
	\begin{equation}
	\| \mathcal{N}_\Omega\|_{op}
	\leq \| \mathcal{N}_B\|_{op} =\frac{4}{\pi^{2}}
	\label{ex3}
	\end{equation}
	for any domain $\Omega$ with $|\Omega|=|B|,$ where  $\Omega^{*}\equiv B$ is the unit ball. Here $\|\cdot\|_{op}$ is the operator norm of the Newton potential on the space $L^2$.
\end{ex}


\begin{thebibliography}{HOHOLT08}

\bibitem{AK}
Arazy J., Khavinson D.: Spectral estimates of {C}auchy's transform in {$L^2(\Omega)$}. Integral Equations Operator Theory, \textbf{15} (6), 901--919, (1992)

\bibitem{AKL}
Anderson J M., Khavinson D., Lomonosov V.: Spectral properties of some integral operators arising in potential
theory. Quart. J. Math. Oxford Ser. (2), \textbf{43}(172), 387--407, (1992)

\bibitem{Ban80}
Bandle C.: Isoperimetric inequalities and applications. Pitman (Advanced Publishing Program), \textbf{7} of {\em Monographs and Studies in Mathematics}, Boston, Mass.-London, (1980)

\bibitem{BLM}
Banuelos R., Latala R., Mendez-Hernandez P.J.: A Brascamp-Lieb-Luttinger-type inequality and applications to symmetric stable processes.
Proc. Amer. Math. Soc., \textbf{129} (10), 2997--3008, (2001)

\bibitem{Ben}
Benguria R.D., Linde H., and Loewe B.: Isoperimetric inequalities for eigenvalues of the {L}aplacian and the
{S}chr{\"o}dinger operator. Bull. Math. Sci., \textbf{2}(1), 1--56, (2012)

\bibitem{BS}
Birman M. {\v{S}}., Solomjak  M. Z.: Estimates for the singular numbers of integral operators. Uspehi Mat. Nauk, \textbf{32} (1(193)), 17--84, (1977)


\bibitem{DR} Delgado J., Ruzhansky M.: Schatten classes on compact manifolds: kernel conditions.
J. Funct. Anal., \textbf{267}(3), 772--798, (2014)

\bibitem{Dos11}
Dostanic M.R.: Regularized trace of the inverse of the Dirichlet Laplacian. Comm. Pure Appl. Math., \textbf{64}(8), 1148--1164, (2011)

\bibitem{Dos12}
Dostanic M. R.: The asymptotic behavior of the singular values of the convolution operators with kernels whose Fourier transform are rational. J. Math. Anal. Appl., \textbf{295}(2), 496--500, (2012)

\bibitem {HH} Harrell E.M., Hermib L.: Differential inequalities for Riesz means and Weyl-type bounds for eigenvalues, J. Funct. Anal., \textbf{254}, 3173--3191, (2008)

\bibitem{GK}
Gohberg I.C.,  Kre{\u\i}n  M.G.: Introduction to the theory of linear nonselfadjoint operators. Translated from the Russian by A. Feinstein. Translations of Mathematical Monographs, \textbf{18}. American Mathematical Society, Providence,
R.I., (1969)

\bibitem{He}
Henrot A.: Extremum problems for eigenvalues of elliptic operators. Frontiers in Mathematics. Birkh{\"a}user Verlag, Basel, (2006)

\bibitem{Kac1}
Kac M.: On some connections between probability theory and differential and
integral equations. In {\em Proceedings of the {S}econd {B}erkeley {S}ymposium on
	{M}athematical {S}tatistics and {P}robability, 1950}, 189--215. University of California Press, Berkeley and Los Angeles, (1951)

\bibitem{Kac2}
Kac M.:
On some probabilistic aspects of classical analysis.
Amer. Math. Monthly, \textbf{77}, 586--597, (1970)

\bibitem{Kac3}
Kac M.:
Integration in function spaces and some of its applications.
Accademia Nazionale dei Lincei, Pisa, (1980).
Lezioni Fermiane. [Fermi Lectures].

\bibitem{Kac4} Kac M.:
An application of probability theory to the study of Laplaces equation.
Ann. Soc. Polon. Math. \textbf{25}, 122--130, (1953)

\bibitem{Kac5} Kac M.:
Distribution of eigenvalues of certain integral operators,
Mich. Math. J., \textbf{3}, 141--148, (1956)

\bibitem{Kac6} Kac M.:
Aspects probabilistes de la Theorie du Potentiel,
Les Presses de la Universite de Montreal, (1970)

\bibitem {KS1} Kal'menov T.Sh., Suragan D.: To spectral problems for the volume potential. Doklady Mathematics, \textbf{80}, 646--649 (2009)

\bibitem{KS2}  Kal'menov T.Sh., Suragan D.: Boundary conditions for the volume potential for the polyharmonic equation. 
Differential Equations, \textbf{48},  604--608, (2012)

\bibitem{KS3} Kalmenov T.Sh.,  Suragan D.: A boundary condition and Spectral Problems for the Newton Potentials.
Operator Theory: Advances and Applications, \textbf{216}, 187--210, (2011)

\bibitem{KS16} Kassymov A., Suragan D.: Some Spectral Geometry Inequalities for Generalized Heat Potential Operators. Complex Analysis and Operator Theory, to appear (2016),
DOI: 10.1007/s11785-016-0605-9

\bibitem{Lan66}
Landkof N.~S.: 
Foundations of modern potential theory.
Springer-Verlag, New York-Heidelberg, (1972).
Translated from the Russian by A. P. Doohovskoy, Die Grundlehren der
mathematischen Wissenschaften, Band 180.

\bibitem{LL}
Lieb E.~H., Loss M.: Analysis, volume~14 of {\em Graduate Studies in Mathematics}.
American Mathematical Society, Providence, RI, second edition, (2001)

\bibitem{Lu}
Luttinger J.~M.: 
Generalized isoperimetric inequalities.
Proc. Nat. Acad. Sci. U.S.A., \textbf{70}, 1005--1006, (1973)

\bibitem{Po}
P{\'o}lya G.:
On the characteristic frequencies of a symmetric membrane.
Math. Z., \textbf{63}, 331--337, (1955)

\bibitem{Po1}
P{\'o}lya G., Szeg{\"o} G.:
Isoperimetric {I}nequalities in {M}athematical {P}hysics.
Annals of Mathematics Studies, \textbf{27}. Princeton University Press,
Princeton, N. J., (1951)

\bibitem{Ruzhansky-Suragan:Newton}
Rozenblum G.,  Ruzhansky M. and Suragan D.:
Isoperimetric inequalities for Schatten norms of Riesz potentials.
J.  Funct. Anal., \textbf{271}, 224--239, (2016)

\bibitem{Ruzhansky-Suragan:PAMS}
Ruzhansky M., Suragan D.:
On Kac's principle of not feeling the boundary for the Kohn Laplacian on the Heisenberg group.
Proc. Amer. Math. Soc., \textbf{144}(2), 709--721, (2016)

\bibitem{Ruzhansky-Suragan:UMN}
Ruzhansky M., Suragan D.:
Schatten's norm for convolution type integral operator.
Russ. Math. Surv., \textbf{71}, 157--158, (2016)

\bibitem{Ruzhansky-Suragan:arxiv}
Ruzhansky M., Suragan D.:
Isoperimetric inequalities for the logarithmic potential operator.
J. Math. Anal. Appl., \textbf{434}, 1676--1689, (2016)

\bibitem{Ruzhansky-Suragan:BMS}
Ruzhansky M., Suragan D.: On first and second eigenvalues of Riesz transforms in spherical and hyperbolic geometries.
Bull. Math. Sci., \textbf{6}, 325-334, (2016)

\bibitem{Ruzhansky-Suragan:Layers}
M.~Ruzhansky and D.~Suragan:
\newblock Layer potentials, {K}ac's problem, and refined {H}ardy inequality on
homogeneous {C}arnot groups.
\newblock {\em Adv. Math.}, 308:483--528, 2017.

\bibitem{Sa}
Saito N.: Data analysis and representation on a general domain using eigenfunctions of {L}aplacian. Appl. Comput. Harmon. Anal., \textbf{25}(1), 68--97, (2008)

\bibitem{Tr67}
Troutman J. L.: The logarithmic potential operator. Illinois J. Math., \textbf{11}, 365--374, (1967)

\bibitem{Tr69}
Troutman J. L.: The logarithmic eigenvalues of plane sets. Illinois J. Math., \textbf{13}, 95--107, (1969)

\bibitem{Zo} Zoalroshd S.: A note on isoperimetric inequalities for logarithmic potentials. J. Math. Anal. Appl., \textbf{437}(2), 1152-1158, (2016)



\end{thebibliography}
\end{document}